\documentclass[10pt, journal, onecolumn]{IEEEtran}

\usepackage[cmex10]{amsmath}
\usepackage{amssymb, amsthm, amsfonts}
\usepackage{graphicx, subfigure} \graphicspath{{figs/}}
\usepackage{fullpage}
\usepackage{enumerate}
\usepackage{comment}

\newcommand{\abcd}{\begin{pmatrix}a&b\\c&d\end{pmatrix}}
\newcommand{\abcdvert}{
\setlength{\arraycolsep}{2.5pt}
\!\begin{pmatrix}
&a&\\[-0.7ex]
b& &c\\[-0.7ex]
&d&\end{pmatrix}}
\DeclareMathOperator{\Tr}{Tr}
\newcommand{\set}[1]{\left\{ #1 \right\} }
\newcommand{\pdi}{$p$-$i$ }
\newcommand{\Ref}[1]{(\ref{#1})}

\newtheorem{theorem}{Theorem}
\newtheorem{lemma}[theorem]{Lemma}
\newtheorem{conjecture}[theorem]{Conjecture}
\theoremstyle{definition}
\newtheorem{defn}{Definition}

\begin{document}
\title{Accurate Lower Bounds on Two-Dimensional Constraint Capacities From
Corner Transfer Matrices}
\author{Yao-ban Chan and Andrew Rechnitzer
\thanks{Manuscript submitted 23/10/2012; revised 6/12/2013.}%
\thanks{YBC is at the University of Queensland (formerly Universi\"at Wien); 
\texttt{yaoban.chan@uq.edu.au}}%
\thanks{AR is at the University of British Columbia; 
\texttt{andrewr@math.ubc.ca}}}

\maketitle

\begin{abstract}
We analyse the capacity of several two-dimensional constraint families --- the
exclusion, colouring, parity and charge model families. Using Baxter's corner
transfer matrix formalism combined with the corner transfer matrix
renormalisation group method of Nishino and Okunishi, we calculate very tight
lower bounds and estimates on the growth rates of these models. Our results
strongly improve previous known lower bounds, and lead to the surprising
conjecture that the capacity of the even and charge(3) constraints are
identical.
\end{abstract}


\section{Introduction}
\IEEEPARstart{I}{n} this paper, we analyse the capacity of several
two-dimensional constraints. A constraint induces a model of magnetic spins,
each of which can take a small number of spin values (typically 2, representing
`up' or `down' magnetic spins). These spins are arranged in an ordered fashion
on a two-dimensional plane --- in all the models we analyse, they are positioned
at the vertices of the square lattice $\mathbb{Z}^2$. Each constraint limits in
some way the potential arrangements of the spins on the lattice.

An application where constraints arise is that of data storage. Here, the values of the spins contain data stored (in some encoding) on a magnetic disk. Due to properties of the disk, it may be undesirable for certain local arrangements of spins to exist, and this leads to a corresponding constraint. For example, in the hard squares model we define below, the spins can take the values $\oplus$ and $\ominus$. However, the presence of two $\ominus$ spins next to each other may cause magnetic interference, so we do not allow this to happen.

Given a constraint such as this, we wish to know what the storage capacity of the disk is. In the language of statistical mechanics, the \emph{partition function} is the number of different configurations that can be written on $N$ spins under the constraint:
\begin{align}
 Z_N = (\mbox{\# of legal configurations on $N$ spins}).
\end{align}
We are interested in the behaviour of this number as the size of the storage
array grows to infinity, and so define the \emph{partition function per site}
\begin{align}
\kappa = \lim_{N \rightarrow \infty} Z_N^{1/N}.
\end{align}
We also call $\kappa$ the \emph{growth rate} of the model. The growth rate is
directly related to the \emph{capacity} of a constraint, which is defined to be
the number of bits of information per spin that is possible to be written on a
disk with that constraint:
\begin{align}
\mbox{cap}(C) = \lim_{N \rightarrow \infty} \frac{1}{N} \log_2 Z_N
= \log_2 \kappa.
\end{align}
In other words, a disk with $N$ spins can store approximately $\mbox{cap}(C) \times N$ bits of information under the constraint $C$. A `free' disk with no constraints has a capacity of 1, since there are $2^N$ possible configurations on $N$ spins, and as such each spin represents one bit of information. As we approach this problem from a combinatorial point of view, we shall limit our terminology to growth rates rather than capacity (of course, they are trivial to derive from each other).

This problem is also studied in symbolic dynamics and dynamical systems, but with different terminology. Here, the capacity is also known as the (topological) \emph{entropy} of the model, while the spins are \emph{letters} (of an alphabet). Constraints are called \emph{sofic shifts} (or systems), with a slight distinction: in constrained coding, the focus is on configurations of finite systems of arbitrary size, while in symbolic dynamics (and this paper), the focus is on configurations on the entire lattice. The exclusion and colouring models below are sub-classes of so-called `shifts of finite type' or `finite memory systems'.

The capacity of constraints in one dimension is, in general, very easily computed using a transfer matrix. Unfortunately, capacities of two-dimensional constraints are significantly more difficult to compute, and because of this we seek estimates and rigorous lower and upper bounds for their values. Of interest here are the methods of Engel \cite{Engel1990}, and Calkin and Wilf \cite{CalkinWilf}, who provide a method to calculate such bounds using transfer matrices, illustrating with the hard squares model. Louidor and Marcus \cite{louidor2010} generalised this method and applied it to the exclusion, parity and charge families of constraints discussed here. Additional references are provided for individual models as we introduce them below.

In this paper, we use corner transfer matrix (CTM) formalism to derive very precise lower bounds on the growth rates of the models that we study. The corner transfer matrix method was first devised by Baxter in 1968 \cite{BaxterDimers} as a way of numerically estimating the growth rate of fully-packed dimers on a rectangular lattice. In later studies in the late 70s \cite{CTM:1, CTM:Hsq}, it was developed into method to calculate series expansions of the partition function and magnetisation. A particularly noteworthy triumph of this method was Baxter's famous solution of the hard hexagons model \cite{Hardhex, Exactly}, which was derived by noticing a pattern in the eigenvalues of the corner transfer matrices.

As far as we know, CTM formalism has never been used before to calculate lower bounds, rather than estimates (although in fact, the theoretical machinery required to produce lower bounds has long been in place). When used for this purpose, the accuracy of the bound depends on how close we can come to the solution of the model.

For this purpose, we use the corner transfer matrix renormalisation group (CTMRG) method of Nishino and Okunishi \cite{CTM:RG1, CTM:RG2}. This method combines the original CTM method with ideas from the density matrix renormalisation group method (see, for example, \cite{White1992, Schollwock2005}), and also produces estimates of the growth rate. The CTMRG is more general and easier to apply numerically than Baxter's original method, and has been applied, among other things, to self-avoiding walks \cite{CTM:RG-FP2} and the Ising model scaling function \cite{Scaling}. Recently, it too has been extended by one of the authors \cite{CTM:Hardsq} to calculate series expansions.

In Section \ref{sec:models}, we define the models that we study and discuss some reformulations needed to make them amenable to CTM analysis. In Section \ref{sec:ctm}, we explain the corner transfer matrix formalism that we use to derive our lower bounds and estimates. These numerical results are presented in Section \ref{sec:results}. In Section \ref{sec:charge-even}, we discuss the (conjectured) equivalency between two of our models, and conclude in Section \ref{sec:conclusion}.

\section{The models}\label{sec:models}

We calculate the growth rate of several models in this paper. All of these models have \emph{spins} which lie on the vertices of a
square lattice, which usually take the values $\oplus$ and $\ominus$.  The
configurations will be constrained in certain ways that we describe below.

\subsection{Exclusion models}
\begin{figure}
 \centering
  \includegraphics[height=5.5cm]{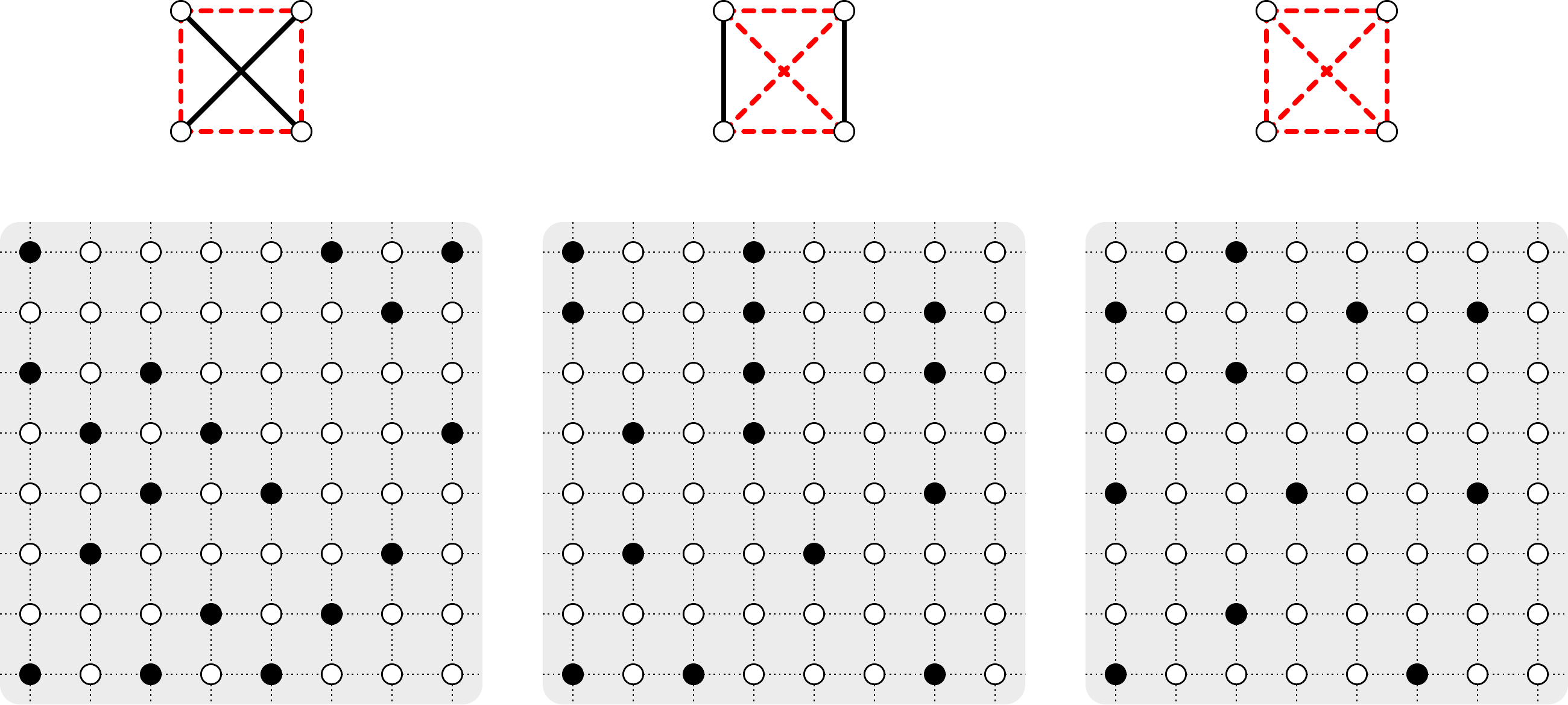}
  \caption{(bottom) Examples of the three exclusion models we consider (from
left to right): hard squares, RWIM and NAK. We represent $\ominus$ spins by 
black vertices and $\oplus$ spins by white vertices. (top) The diagrams
show the allowed configurations around a face for each model; two vertices
joined by a solid line may be simultaneously occupied by $\ominus$s, while vertices joined by
a dashed line may not.}
\label{fig exclusion eg}
\end{figure}
In the exclusion models, we consider configurations of $\oplus, \ominus$ spins which
are forbidden to contain certain local configurations of $\ominus$s (see
Fig.~\ref{fig exclusion eg}).
\begin{defn}
Consider a configuration of $\oplus,\ominus$ spins on the vertices of the
square lattice.
\begin{itemize}
\item The configuration satisfies the \emph{hard squares} constraint if it does
not contain two $\ominus$ spins joined by a single bond. If one considers the
$\ominus$ spins to be the centres of rotated square blocks, then this
constraint requires that these blocks do not overlap. 

This model is a special case of a more general hard squares model where the constraint is enforced, but each $\ominus$ spin is also given a weight $z$. Our model corresponds to $z = 1$. In the more general model, questions of interest include not only this special case, but also the location of the phase transitions (singularities in the free energy). These questions have been studied using corner transfer matrices in \cite{CTM:Hsq, CTM:Hardsq}. Placing the spins on the vertices of a triangular lattice results in the hard hexagon model, which has been exactly solved using corner transfer matrices \cite{Hardhex, Exactly}. Studies on the $z=1$ case include \cite{CalkinWilf, BaxterPlanar}.

\item The configuration satisfies the \emph{read-write isolated memory} (RWIM)
constraint \cite{RWIM1, RWIM2} if it does not contain two $\ominus$ spins joined by a horizontal
bond, or lying diagonally across a single face. This constraint can be used to model one-dimensional memory states with the hard square constraint, which are sequentially updated in such a way that no two adjacent spins are changed in the same update.

\item The configuration satisfies the \emph{non-attacking kings} (NAK)
constraint if it does not contain two $\ominus$ spins joined by a single bond,
or lying diagonally across a single face. Equivalently every face of the
lattice contains at most one $\ominus$ spin. If one considers each $\ominus$
spin to be a king on a chessboard, then the NAK constraint requires that no
king can attack another. This has been studied, along with the hard squares and other models, in \cite{Weeks}.
\end{itemize}

\end{defn}

\subsection{Colouring models}

We also consider vertex colourings of the lattice; each vertex is coloured using labels from 0 to $q-1$ for some fixed $q$, and again we prohibit certain local configurations.
\begin{defn}
A configuration of $\{0,\ldots,q-1\}$ spins satisfies the \emph{$q$-colouring} constraint if it does not contain two nearest-neighbour spins with the same colour.
\end{defn}

This is related to the chromatic polynomial and the anti-ferromagnetic Potts
model (see for example \cite{Biggs, Shrock, Lundow, Jacobsen}). When $q=2$, there are only two possible configurations,
while for $q=3$ it is known that the growth rate is exactly $(4/3)^{3/2} =
1.539\,600\,7178\dots$ \cite{Lieb1967}. We study this model for higher $q$, where
the growth rate is not known exactly.

\subsection{Parity models}
Again we place the spins $\oplus$ and $\ominus$ at the vertices of the square
grid. However unlike the models above, the constraints are no longer local. The
two parity models arise by constraining the parity of the number of $\oplus$
spins lying between two $\ominus$ spins along any horizontal or vertical line. 
\begin{defn}
Consider a configuration of $\oplus, \ominus$ spins on the vertices of the
square grid. Take any two $\ominus$ spins lying along a horizontal
or vertical line so that there is no other $\ominus$ spin between
them.
\begin{itemize}
 \item If every such pair is separated by an even number of $\oplus$ spins, then
the configuration satisfies the \emph{even} constraint.
\item Conversely, if every such pair is separated by an odd number of $\oplus$
spins, then the configuration satisfies the \emph{odd} constraint.
\end{itemize}
\end{defn}
While the above spin constraints are not local, it is possible to map them to
other models of states with local constraints, which makes them amenable to CTM
analysis. We explain this mapping below.

The growth rate of the odd model is known to be $\sqrt{2}$
\cite{louidor2010}. To prove that the growth rate is at least $\sqrt{2}$,
populate the lattice with $\oplus$ in a checkerboard pattern. The remaining
spins, which constitute half the total spins, can then be assigned $\oplus$ or $\ominus$ independently. This always produces a valid configuration, so the bound follows.

To obtain the opposite bound, consider configurations of spins which obey the
odd constraint along horizontal lines, while having no constraint along
vertical lines. Being less constrained, the growth rate of this new model is
greater than that of the odd model and is equal to the growth rate of the
odd model in 1 dimension. It is straightforward to show that this
is exactly $\sqrt{2}$.

Since the growth rate of the odd model is known exactly, we do not discuss it
further. 

\subsection{Charge(3) model}

A final family of constraints is parameterised by a positive integer
corresponding to the maximum cumulative electric `charge' along a horizontal or
vertical line. As we scan along a row or column, our cumulative `charge'
will increase by $+1$ or $-1$ each time we pass a $\oplus$ or $\ominus$ 
(respectively). We arrive at the charge model by constraining the 
allowed values of this cumulative charge.
\begin{defn}
 A configuration of $\oplus, \ominus$ spins satisfies the \emph{charge(3)}
constraint if for any finite horizontal or vertical line segment, the sum
of the spins along that segment lies between $-3$ and $+3$. 

Equivalently, starting from any given horizontal (resp. vertical) bond, there 
exists an initial charge between 0 and 3, so that scanning to the east 
(resp. south), the cumulative charge always lies between $0$ and $3$.
\end{defn}

Clearly the above definition can be generalised by replacing $3$ by any other
positive integer. There are only 2 valid charge(1) configurations --- checkerboard patterns of $\oplus$ and $\ominus$ --- so this system has growth rate 1.
The growth rate of charge(2) configurations on $\mathbb{Z}^d$ was proved in
\cite{louidor2010} to be $2^{2^{-d}}$. Charge(3) is the first of this family of
models for which the growth rate is unknown --- indeed, previous works have been
unable to determine its first digit exactly. Our results strongly support the
following surprising conjecture:
\begin{conjecture}
\label{conj charge even}
 The charge(3) and even models have equal growth rate on $\mathbb{Z}^2$.
\end{conjecture}

\subsection{Model transformations: states on bonds}

The exclusion and colouring models all share the property that their constraint 
is local, acting entirely within the confines of a single cell. This is 
precisely what we need for the model to be amenable to analysis using corner 
transfer matrices. As such, we can apply this analysis to these models directly.

However, the even and charge(3) constraints are non-local. Fortunately, we can 
transform these models so that the constraints can be expressed in a local 
manner. To do this, we will map configurations of $\oplus$ and $\ominus$ spins 
on the vertices of the lattice to configurations of \emph{states} that lie on 
the bonds of the lattice.

We first explain this process for the charge(3) model (as it is conceptually simpler). Here, the state of any horizontal (resp. vertical) bond 
is the cumulative charge on that bond, having read the spin to its west (resp. north). So in any 
valid configuration, all bonds will have integer states between 0 and 3 (inclusive). Considering the states around a vertex spin, the states of the south and east bonds must be 1 more or 1 less than those of the north and west bonds respectively (and the differences must be consistent with each other). That is,
\begin{align}
 \begin{array}{ccc}
  &a&\\
  b\!\!&\circledast&\!\!c \\
  &d& 
 \end{array}
  & = \begin{cases}
      \text{valid } & \text{ if } d-a=c-b=\pm 1\\
      \text{invalid } & \text{ otherwise. }
    \end{cases}\label{eq:charge_weight}
\end{align}
We call this bond model the $Q$-charge model ($Q$ being traditional
notation for charge in electrostatics).

Note that this mapping is not a bijection; it is possible for a single charge(3) configuration to correspond 
to many $Q$-charge configurations. If the cumulative charges along a given row 
take 4 different values, then only a single $Q$-charge configuration is 
possible on that row. However if the cumulative charges take only 2 or 3 values, then it is 
possible to increase or decrease those values and so create 
several possible $Q$-charge configurations along the row. However, these 
possibilities are not sufficient to change the growth rate of the model.

\begin{lemma}
  The growth rate of the charge(3) and $Q$-charge models are equal.
\end{lemma}
\begin{IEEEproof}
  The mapping defined above is one-to-many and so shows that $\kappa_{charge(3)} \leq
  \kappa_{\text{$Q$-charge}}$. To prove the reverse inequality, consider both
constraints on the lattice $[1,N]^2$ (including boundary bonds).

Starting with a valid charge(3) configuration, the state of the first bond in 
each row/column may have up to 4 possible values (in fact, one of these values is always impossible, but this does not have an effect on the proof). After this, the row/column is fully determined. So on this 
finite lattice, every charge(3) configuration corresponds to at most $4^{2N}$ 
$Q$-charge configurations. Since there are $N^2$ vertices in this lattice, the 
growth rates are equal in
the $N \to \infty$ limit.
\end{IEEEproof}

Now we apply a similar conversion to the even model. Here, the state of a given horizontal (resp. vertical) 
bond records the number of $\oplus$ spins lying between that bond and the first 
$\ominus$ to its west (resp. north). However, we do not need to know the total 
number of $\oplus$ passed, but only whether it is an even or odd number. 
Hence the states of each bonds will be `$p$' or `$i$' taken from the French 
\emph{pair} and \emph{impair} for even and odd. We chose this convention to 
avoid confusion between the names of these bond states and the names of the 
even and odd models. 

All bonds adjacent to $\ominus$ spins are labelled `$p$'. Around any $\oplus$ 
spin, bonds on opposite sides of the vertex must take different
labels --- so a bond with a `$p$' state must lie opposite a bond with an
`$i$' state and vice versa. This means that for any given $\ominus$ spin, if we 
scan east to the next $\ominus$ spin, we will pass an even number of $\oplus$
spins and so read an odd-length alternating sequence of bond states $pipi\dots
ip$. 

We call this new model the `$p$-$i$' model. A labelling of bonds by $p$ and $i$ 
will be valid for this model when there are no two adjacent $i$ states along any 
horizontal or vertical line, and any vertex must be adjacent to zero or two $i$ 
states.

Again, this mapping is not bijective, as there are two possible \pdi 
configurations for a row which contains only~$\oplus$~spins in the even model, 
namely $pipi\dots$ or $ipip\dots$.  However, it is again easy to show that the 
growth rate is unchanged.

\begin{lemma}
  The growth rate of the even model is equal to that of the \pdi model.
\end{lemma}
\begin{IEEEproof}
  The proof is identical to that of the previous lemma, except that there are 
at most two possible ways to determine each row or column. This means that 
every even configuration corresponds to at most $2^{2N}$ \pdi configurations, 
which is still sub-exponential in the number of vertices in the lattice 
(which is $N^2$). 
\end{IEEEproof}


\begin{figure}
 \centering
  \includegraphics[width=0.8\textwidth]{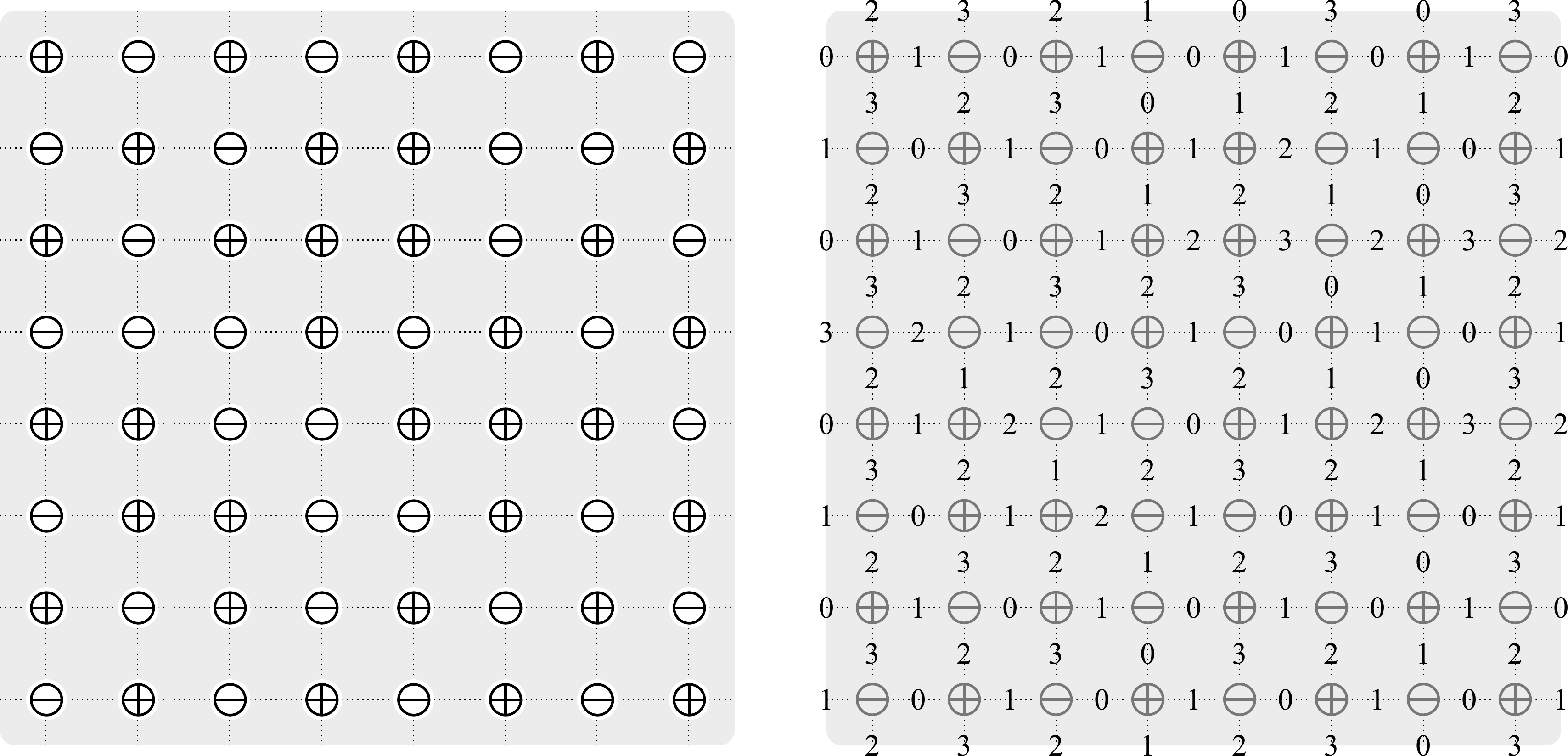}
  \caption{An example of a valid charge(3) configuration (left) and one of the 
$Q$-charge configurations that it maps to (right).}
\end{figure}
\begin{figure}
 \centering
  \includegraphics[width=0.8\textwidth]{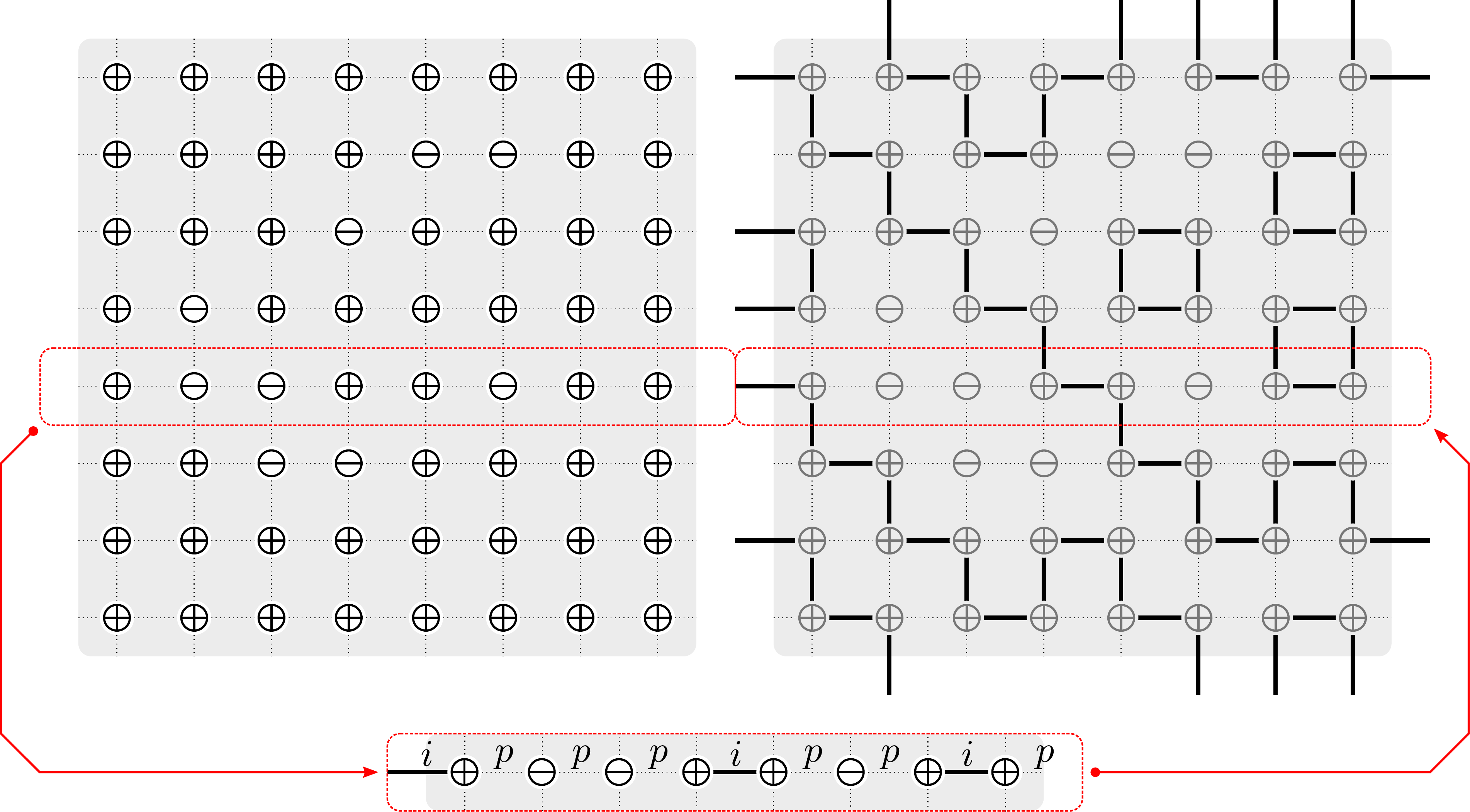}
  \caption{An example of a valid even configuration (left) and one of the \pdi 
configurations it maps to (right). In the under-set figure we show 
explicitly how one row of the even configuration is mapped to a row in the \pdi 
model. Bonds between $\ominus$ spins are labelled by an odd length 
sequence $pipi\dots\ p$. To visualise this we occupy each bond labelled by $i$. 
Notice that each vertex has degree 0 or 2, and that if two occupied bonds meet 
at a vertex, then they must be orthogonal to each other. }
\label{fig pi eg}
\end{figure}

\subsection{Model transformations: states on faces}\label{sec even vert}

We make a further transformation from the \pdi model to another model (which we 
call \emph{even-face}), which places states on the centres of the faces of the 
lattice. These states can take the values 0 or 1. It is not necessary to make 
this transformation, but empirically it responds better to the corner transfer 
matrix renormalisation group method.

Consider each bond in state `$i$' in a \pdi configuration to be occupied 
(see Fig.~\ref{fig pi eg}). Each vertex in a valid configuration is 
incident to either 0 or 2 occupied bonds, and these form 
distinct simple loops which separate the plane into regions. We now properly colour these regions with 2 colours, and define the state of each face to be its colour.

\begin{figure}\centering
  \includegraphics[width=0.9\columnwidth]{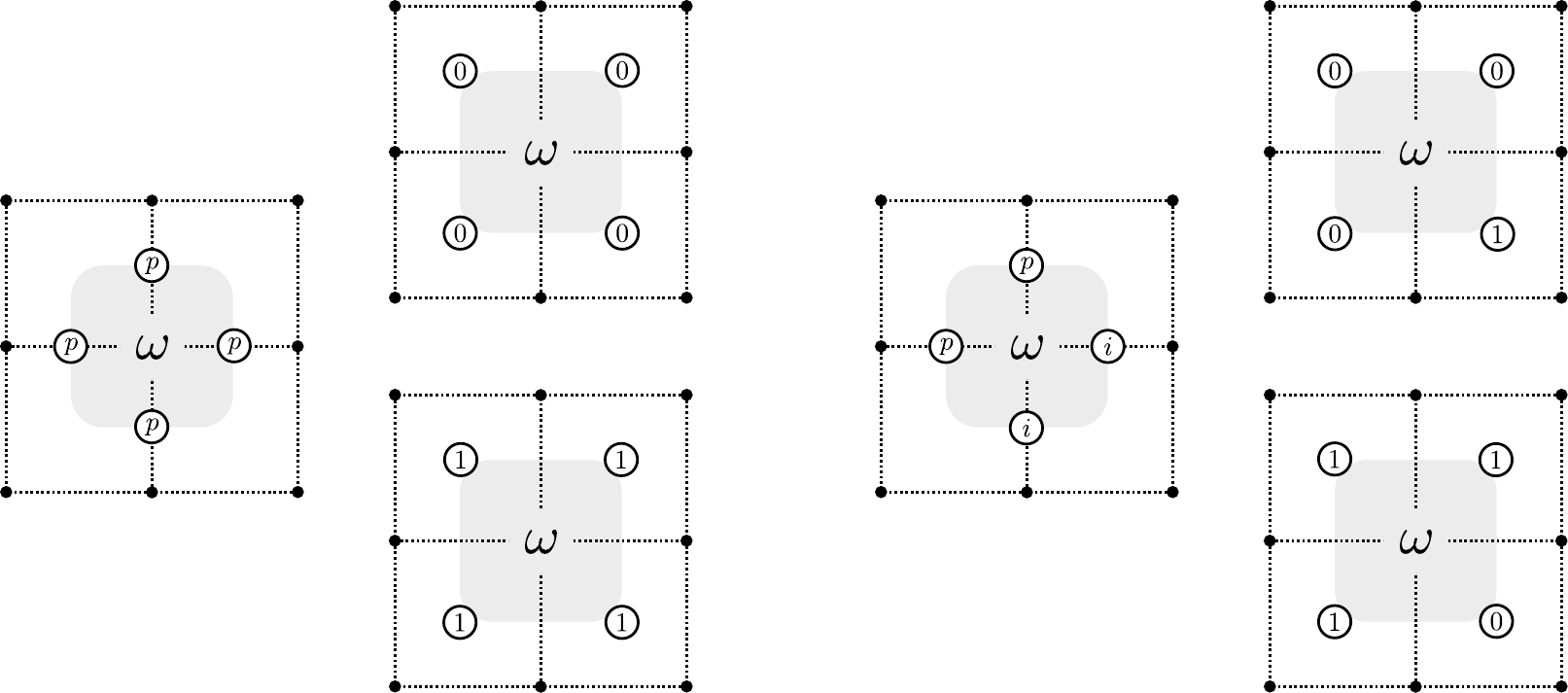}
\caption{Possible cells in the \pdi and even-face models, up to rotation.}\label{fig:even_b2f}
\end{figure}

Thus, given a \pdi configuration, two states on neighbouring faces in the 
corresponding even-face configuration are equal if they are separated by a $p$ 
bond (since they are part of the same region), and unequal otherwise (since 
they are separated by an occupied bond and so in different regions). A full 
enumeration of the possible \pdi configurations shows that the resulting 
configuration is internally consistent (see Fig. \ref{fig:even_b2f}). Now we 
determine the validity of a configuration according to each cell:
\begin{align}
\abcd &= \begin{cases} 
\mbox{invalid} & \mbox{if $a+b+c+d=2$} \\
\mbox{valid} & \mbox{otherwise.} 
\end{cases}
\end{align}

Again, this mapping is not a bijection, but determining the value of one state on an arbitrary face in the even-face model will fully determine the rest of the configuration. As such, the number of even-face configurations is exactly twice the number of \pdi configurations, and their growth rates are equal.

\section{Corner transfer matrix formalism}\label{sec:ctm}

 \subsection{Computing a rigorous lower bound}

  Consider, for the sake of exposition, a model in which spins lie on the
vertices of $\mathbb{Z}^2$ and can take two values. Further, the model 
is governed by an `interaction round a face' $\omega$ that is symmetric across a vertical line. This includes the exclusion models (hard
squares, NAK and RWIM), the colouring models and the even-face model (since the
square lattice is self-dual). More work is required for the charge model and we discuss it further below. We wish to count the number of configurations, and so take
\begin{align}
  \omega\!
\begin{pmatrix}
  a&b\\c&d
\end{pmatrix}
&= \begin{cases}
    1 & \text{if the local configuration is valid}\\
    0 & \text{otherwise.}
   \end{cases}
\end{align}

\begin{figure}
\centering
 \includegraphics[height=5cm]{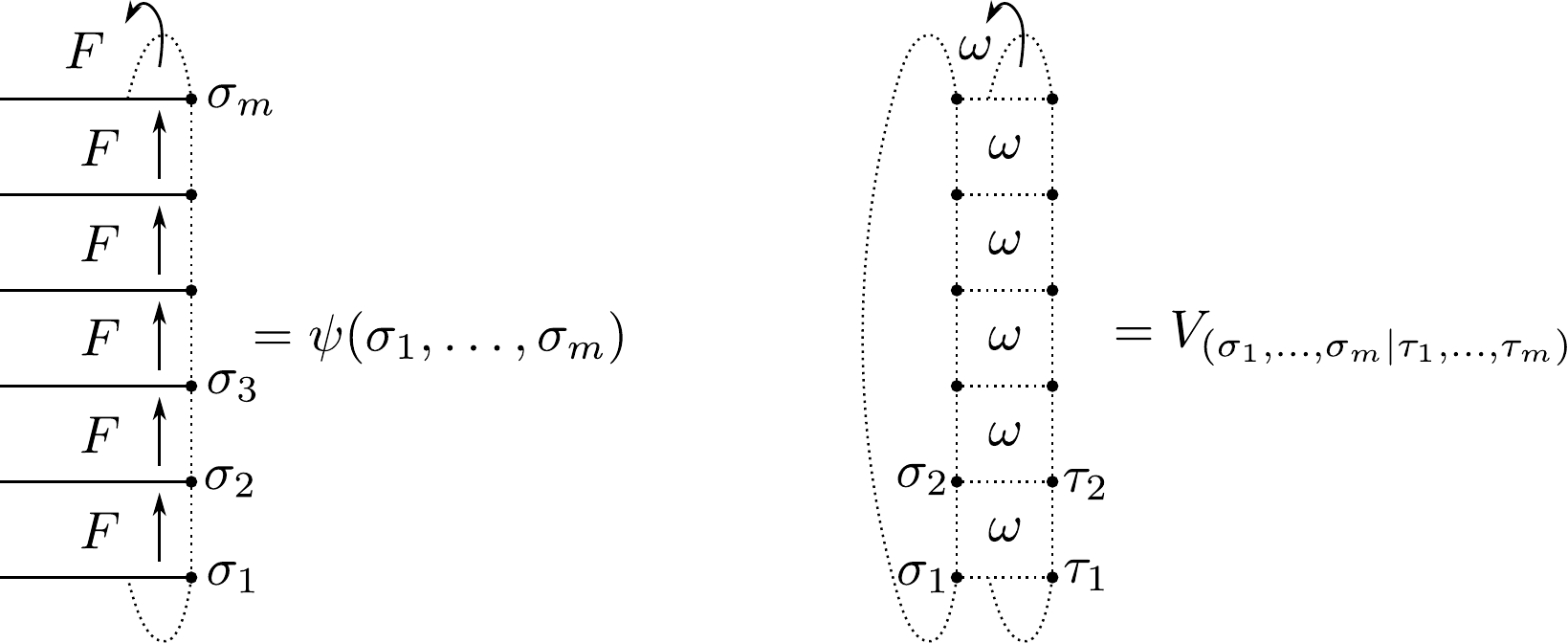}
 \caption{We define the $2^m$ elements of the vector $\psi$ to be the trace of
the product of $F$ matrices. The elements of the column transfer matrix $V$
are products of the face weight $\omega$. In both cases we apply cylindrical
boundary conditions so that $\sigma_{m+1} \equiv \sigma_1$ and $\tau_{m+1}
\equiv \tau_1$.}
  \label{fig vandf}
\end{figure}

First consider a strip of the square lattice of height $m$ with
cylindrical boundary conditions identifying the top and bottom edges. Let
$V$ be the $2^m \times 2^m$ column transfer matrix defined by
\begin{align}
 V_{\sigma, \tau} &= \prod_{i=1}^m \omega \!\begin{pmatrix}
                                         \sigma_{i+1} & \tau_{i+1}\\
                                         \sigma_i & \tau_i\\
                                        \end{pmatrix}
\end{align}
where $\sigma_{m+1} \equiv \sigma_1$ and similarly for $\tau$. See
Fig.~\ref{fig vandf}. Note that this is just the traditional column transfer matrix
and is symmetric.

Since $V$ is symmetric, we can apply a variational principle --- namely, the
dominant eigenvalue $\Lambda(m)$ is given by
\begin{align}
 \Lambda(m) &= \max_{\psi} \frac{\psi^T V \psi }{\psi^T \psi}.\label{eq:var}
\end{align}
The growth rate of the system is the limit
\begin{align}
  \kappa &= \lim_{m \to \infty} \Lambda(m)^{1/m}.
\end{align}
To obtain a lower bound, it suffices to substitute a particular $\psi$ into (\ref{eq:var}), which we do
using Baxter's corner transfer matrix ansatz. 

Let $\set{ F(a,b) }$ be a set of four $n\times n$ matrices 
indexed by two spin values $a,b \in \set{0,1}$ so that $F(a,b)= F^T(b,a)$.
While some of these matrices may be zero (e.g. for the hard squares model
$F(1,1) = \mathbf{0}$), they should not all be trivial. Consider the vector of
dimension $2^m$ which satisfies
\begin{align}
  \psi_{\sigma} &= \Tr \left[
  F(\sigma_1, \sigma_2) F(\sigma_2, \sigma_3) \dots F(\sigma_m, \sigma_1)
  \right].
\end{align}
This is shown in Fig.~\ref{fig vandf}. For any particular choice of $F$, we
obtain a vector $\psi$ and so a lower bound for the maximisation problem. Now,
given a choice of $F$, we need to be able to compute $\psi^T V \psi$ and $\psi^T
\psi$ and so our lower bound.

\begin{figure}
\centering
 \includegraphics[height=8cm]{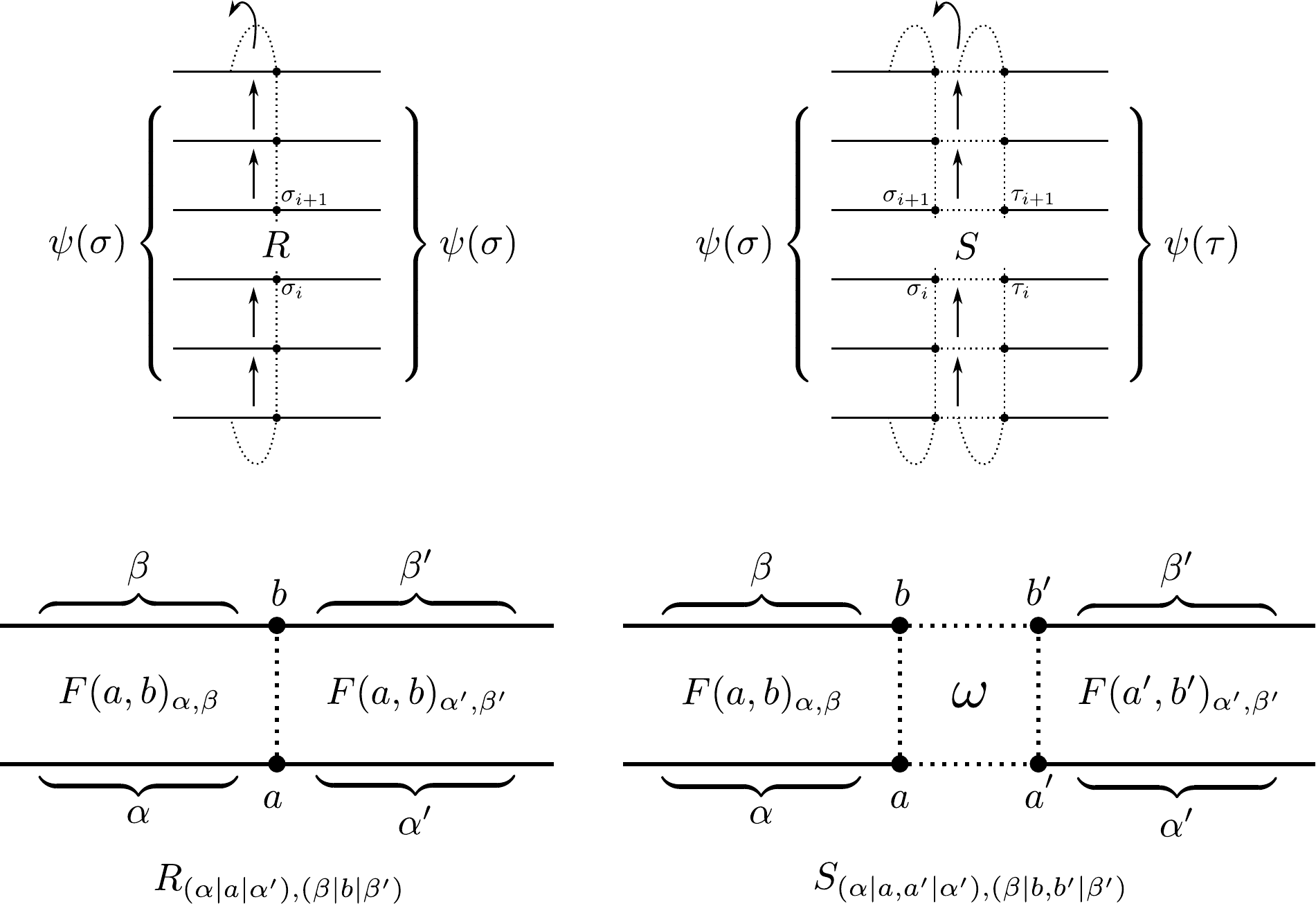}
\caption{In order to compute $\psi^T\psi$ (top-left) and $\psi^T V \psi$
(top-right), we define two matrices $R$ and $S$ so that $\psi^T \psi = \Tr(R^m)$
and $\psi^T V \psi = \Tr(S^m)$. The elements of $R$ and $S$ are defined as
indicated (bottom-left and bottom-right respectively).}
  \label{fig rands}
\end{figure}

Let us fix $m$ and the matrices $F$. We denote
the $n^2$ entries of $F(a,b)$ as $F(a,b)_{\alpha,\beta}$ for $\alpha,\beta \in
\set{1,\dots,n}$. It is possible to rewrite the product of $\psi^T \psi$ as the trace of a power of a new matrix, $R$. This is extremely helpful because $R$ is not dependent on the height of the strip. Instead, the height dependence is confined to the power of the matrix. In particular, we may write
\begin{align}
 \psi^T \psi &= \Tr( R^m )
\end{align}
where $R$ is a $(2n^2) \times (2n^2)$ matrix defined by
\begin{align}
 R_{(\alpha|a|\alpha'), (\beta|b|\beta')} 
&= F(a,b)_{\alpha,\beta} F(a,b)_{\alpha',\beta'}.
\end{align}
See Fig.~\ref{fig rands} for an illustration of this. To see this consider

\begin{subequations}
\begin{align}
  \psi^T \psi 
  &= \sum_{\sigma} \left( \Tr \prod_{i=1}^m F(\sigma_{i},\sigma_{i+1})
\right)^2\\
  &= \sum_{\sigma} \left( 
\sum_{\alpha_1,\dots,\alpha_m} \prod_{i=1}^m F(\sigma_i,
\sigma_{i+1})_{\alpha_i, \alpha_{i+1}} \right)^2\\
  &= \sum_{\sigma} 
\sum_{\alpha_1,\dots,\alpha_m}
\sum_{\alpha_1',\dots,\alpha_m'}
 \prod_{i=1}^m  F(\sigma_i,\sigma_{i+1})_{\alpha_i,
\alpha_{i+1}}
 F(\sigma_i,\sigma_{i+1})_{\alpha_i', \alpha_{i+1}'}.
\intertext{Now rewrite the summands as elements of $R$}
  &= 
\sum_{\sigma} 
\sum_{\alpha_1,\dots,\alpha_m}
\sum_{\alpha_1',\dots,\alpha_m'}
\prod_{i=1}^m R_{(\alpha_i|\sigma_i|\alpha_i'),
(\alpha_{i+1}|\sigma_{i+1}|\alpha_{i+1}')}\\
&= \Tr( R^m ).
\end{align}
\end{subequations}

In a similar way,
\begin{align}
 \psi^T V \psi &= \Tr(S^m)
\end{align}
where $S$ is a $(4n^2) \times (4n^2)$ matrix defined by
\begin{align}
 S_{(\alpha|a,a'|\alpha'), (\beta|b,b'|\beta')} & =  
F(a,b)_{\alpha,\beta}
\cdot
\omega\begin{pmatrix} b&b'\\a&a' \end{pmatrix} \cdot
F(a',b')_{\alpha',\beta'}.
\end{align}
Note that due to the symmetry $F(a,b) = F^T(b,a)$, both $R$ and $S$ are
symmetric.

Now let $\xi$ and $\eta$ be the maximal eigenvalues of $R$ and $S$ respectively. Our lower
bound becomes
\begin{align}
 \Lambda(m) & \geq \frac{\Tr( S^m ) }{\Tr(R^m)} & \Rightarrow
  \kappa & = \lim_{m \to \infty} \Lambda(m)^{1/m} \geq
\frac{\eta}{\xi}.\label{kappa bound}
\end{align}
So to compute our bound it suffices to compute the dominant eigenvalues of $R$
and $S$.

\begin{figure}
\centering
 \includegraphics[height=4.5cm]{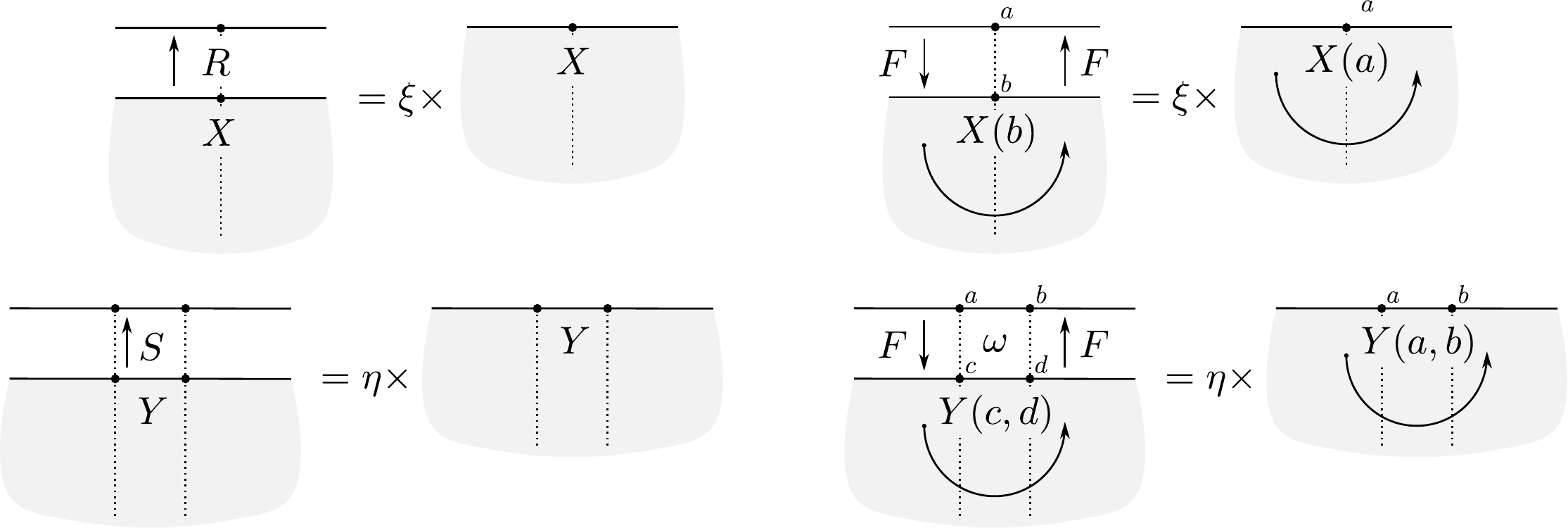}
\caption{(left) $X$ and $Y$ are the dominant eigenvectors of the matrices $R$ and
$S$ (respectively) with eigenvalues $\xi$ and $\eta$ (respectively). (right) By
recasting $X$ and $Y$ as matrices, we can rewrite the eigenvalue equations as
matrix equations in terms of $F$ as shown.}
\label{fig xandy}
\end{figure}

 Let $X$ be the dominant eigenvector of $R$, so that
\begin{align}
 \xi X &= R X.
\end{align}
See Fig.~\ref{fig xandy}. It is helpful in what follows to recast $X$ as a set of matrices, and the
corresponding eigenvalue equation as a matrix equation. The vector $X$ has
dimension $2n^2$ and we index its elements as $X_{(\alpha|a|\beta)}$. Define $X(a)_{\alpha,\beta} = X_{(\alpha|a|\beta)}$. Then we can write
the eigenvalue equation as
\begin{subequations}
\begin{align}
 \xi X(a)_{\alpha,\alpha'} &= (R X)_{(\alpha|a|\alpha')}\\
  &= \sum_{\beta,b,\beta'} R_{(\alpha|a|\alpha'),(\beta|b|\beta')}
  X_{(\beta|b|\beta')}\\
  &= \sum_{\beta,b,\beta'} F(a,b)_{\alpha,\beta} F(a,b)_{\alpha',\beta'}
  X(b)_{\beta,\beta'}\\
  &= \sum_{\beta,b,\beta'} F(a,b)_{\alpha,\beta}
 \,X(b)_{\beta,\beta'}\, F(b,a)_{\beta',\alpha'}  \\
  &= \left(\sum_b F(a,b) \,X(b)\, F(b,a)\right)_{\alpha,\alpha'}
\end{align}
\end{subequations}
and so
\begin{align}
  \xi X(a) &= \sum_b F(a,b) \, X(b) \, F(b,a).
\label{eq:ctm1}
\end{align}
Now define $Y$ to be the dominant eigenvector of $S$, so $\eta Y = SY$. It is of
dimension $4n^2$, and we index its elements as $Y_{(\alpha|a,b|\beta)}$. Again, partition
these as $Y(a,b)_{\alpha,\beta} = Y_{(\alpha|a,b|\beta)}$. With similar working, we have the
eigenvalue equation
\begin{align}
 \eta Y(a,b) &= \sum_{c,d} \omega\begin{pmatrix} a&b\\c&d \end{pmatrix}
  F(a,c) \,Y(c,d)\, F(d,b).\label{eq:ctm2}
\end{align}

To compute $\xi$ and $\eta$, one could apply the power method to the matrices
$R$ and $S$ using a random initial vector. However, since $R$ and $S$ have
dimension $O(n^2)$, each matrix-vector multiplication (and so each
iteration) would take $O(n^4)$ operations. Instead it is more efficient to apply
the power method implicitly using \Ref{eq:ctm1} and \Ref{eq:ctm2}.

Consider equation~\Ref{eq:ctm1}. Given a
prospective eigenvector $X_i(a)$ with a fixed normalisation, we substitute this
into the right-hand side of (\ref{eq:ctm1}). Given this normalisation, the next
iteration $\xi_{i+1}, X_{i+1}(a)$ is now given by the left-hand side of the
equation. Each iteration of these equations involves multiplying square matrices
of dimension $n$, and so requires $O(n^3)$ operations (or smaller with the use of more sophisticated algorithms).

We can use any value of $X(a)$ to start the power method, but the CTMRG method
below will suggest an appropriate starting value. To cover the possibility that
this starting value is a sub-dominant eigenvector, we perturb it slightly before
starting the power method. $\eta$ (and $Y$) is found in an identical manner,
except this perturbation is not needed (the result will still give a lower bound
on $\eta$ and so $\kappa$).

So our process for computing a rigorous lower bound for $\kappa$ is this: we
take any set of $F$ matrices, apply the power method to \Ref{eq:ctm1}
and~\Ref{eq:ctm2} to calculate $\eta$ and $\xi$, and then our lower bound is
$\eta/\xi$. 

Of course, the efficacy of this approach depends critically on the set of $F$
matrices --- thus we now need a heuristic to find a `good' set of $F$
matrices. This is achieved by use of the corner transfer matrix renormalisation
group algorithm, which we detail in the next section.

We note that the method described in Sections 3 and 6 of \cite{louidor2010} is similar in spirit to that described above, in that they both replace the eigenvector of the transfer matrix with an expression that simplifies calculations.

\subsection{Corner transfer matrix renormalisation group method}
\label{sec ctmrg}
We now wish to choose a set of $F$ matrices which provide a good lower bound, i.e. as large as possible. We do this using a heuristic algorithm
known as the corner transfer matrix renormalisation group (CTMRG) method. Note
that while the choice of $F$ alters the quality of the bound, it does
not affect the fact that it is a lower
bound, and as such it is unnecessary to
prove any convergence properties of the CTMRG method.

\begin{figure}
\centering
 \includegraphics[width=\columnwidth]{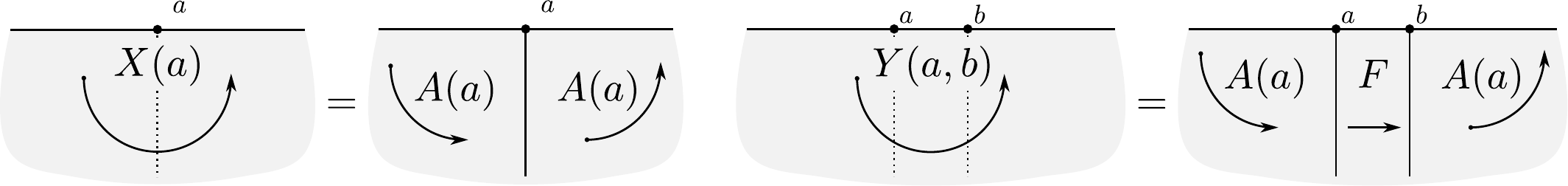}
\caption{In Fig.~\ref{fig xandy} the eigenvectors $X$ and $Y$ were recast as
matrices. Now we rewrite them again as products of the matrix $F$ together with
the corner transfer matrix $A$. The matrix $F$ sweeps through a half-row or
half-column and the matrix $A$ sweeps through a quarter of the lattice.}
\label{fig xytoa}
\end{figure}

Since we wish our lower bound to be as large as possible, we wish to maximise
$\eta/\xi$ with respect to $F$. The derivation of equations which ensure this is
long and complicated, and we shall not repeat it here. It is sufficient \cite{CTM:1}
 to note that if we have a set of $n \times n$ matrices $A(a)$ which satisfy
\begin{align}
 X(a) &= A^2(a) & Y(a,b) &= A(a) \,F(a,b)\, A(b)\label{eq:ctm3}
\end{align}
then $\kappa$ is stationary with respect to $F$. Note that defining $A$ in this
way assumes that the model is invariant under rotation by $\pi/2$. For models
without this symmetry, we define several different corner transfer matrices ---
see the next section. Unfortunately \Ref{eq:ctm3} does not
prove that $\kappa$ is maximal, merely stationary, though empirically this does appear to be the
case. See Fig.~\ref{fig xytoa} for a graphical interpretation of the above
equations. Unlike the more familiar column and row transfer matrices, the matrix
$A(a)$ (the eponymous corner transfer matrix) adds the weight of one quadrant or corner of the lattice.

The CTMRG method calculates matrices which converge to the finite-size solution
of (\ref{eq:ctm1}), (\ref{eq:ctm2}) and (\ref{eq:ctm3}), known collectively as
the \emph{CTM equations}. As such, the $F$ matrices that it produces can be used
in our calculations above to provide (what appears to be) a very tight lower bound.

To gain an intuition for the CTMRG method, consider $A(a)$ and $F(a,b)$ to be
infinite-size transfer matrices (with an appropriate normalisation) which
represent a corner and half-row of the entire plane respectively. It is easily seen
that with $X$ and $Y$ determined according to (\ref{eq:ctm3}), they represent
(normalised) half-plane transfer matrices. In this case,
(\ref{eq:ctm1})~and~(\ref{eq:ctm2}) are seen to be satisfied by this choice of
$F$, $X$ and $Y$. So the infinite-dimensional solution of the CTM equations can
be thought of as a set of transfer matrices.

\begin{figure}\centering
  \includegraphics[height=3cm]{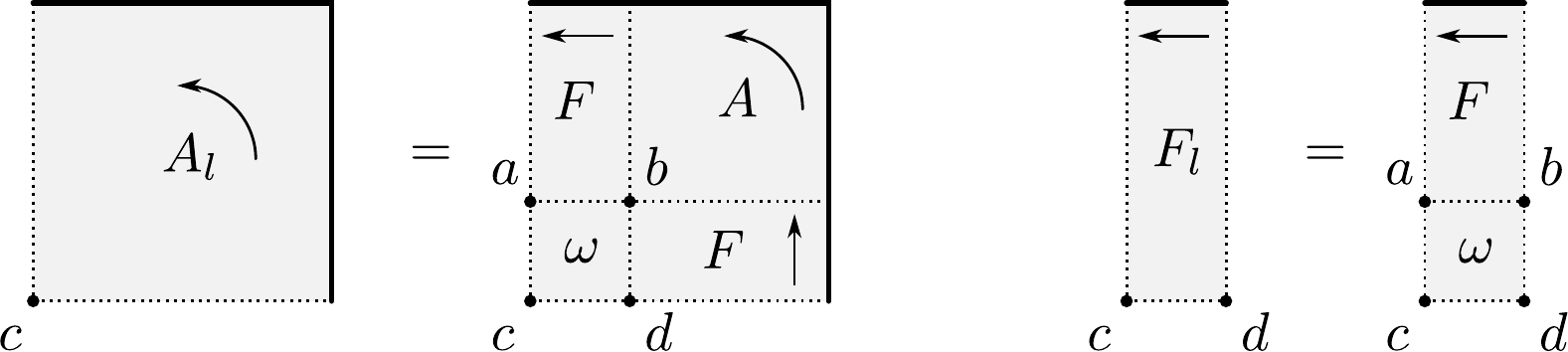}
\caption{Expansion of transfer matrices.}\label{fig:expand}
\end{figure}

In order to approximate this infinite-dimensional solution, we can construct
initial $A$ and $F$ matrices as the corner and half-row transfer matrices of a 
lattice of $(2p+1) \times (2p+1)$ spins. These matrices have dimension $2^p
\times 2^p$.
We can now bootstrap to a larger lattice by constructing the block matrices
\begin{subequations}
\label{eq:expand}
\begin{align}
 A_l(c)|_{d,a} &= \sum_{b} \omega \abcd F(d,b)\, A(b)\, F(b,a) \\ 
 F_l(d,c)|_{b,a} &= \omega \abcd F(b,a)
\end{align}
\end{subequations}
each having dimension $2^{p+1} \times 2^{p+1}$. From Fig.~\ref{fig:expand}, it
is apparent
that $A_l$ and $F_l$ are the corner and half-row transfer matrices of a lattice
of $(2p+3)\times(2p+3)$ spins. By repeating this process indefinitely, we
approach the infinite-size transfer matrices which we require.

At any $p$, we can compute the number of configurations on the
lattice of size $(2p+1)\times(2p+1)$ as $Z_{2p+1,2p+1} = \sum_a \Tr\left( A^4(a)
\right)$. Unfortunately, since the size of these matrices grows exponentially,
this is not a practical approach for even moderate $p$. It does,
however, suggest that we should examine the eigenvalues of $A$. If the
eigenspectrum of $A$ decays quickly, then we can approximate this trace by using
only the largest few eigenvalues of $A$.

To manipulate the eigenvalues directly, we diagonalise the matrix $A_l$. Observe
that the CTM equations are invariant under the similarity transforms
\begin{subequations}
\label{eq:reduce}
\begin{align}
  A(a) &\mapsto P^T(a) \,A(a)\, P(a)
  &F(a,b)& \mapsto P^T(a) \,F(a,b)\, P(b)   \\
  X(a) &\mapsto P^T(a) \,X(a)\,P(a)
  &Y(a,b) &\mapsto P^T(a) \,Y(a,b)\, P(b)
\end{align}
\end{subequations}
for orthogonal matrices $P(a)$. This means that we can take $A_l(a)$ to be
diagonal, with entries ordered from
largest to smallest. We throw away the smaller eigenvalues by truncating 
our large matrices to size $n \times n$. Since the discarded eigenvalues are small,
the effect of this truncation on the trace (and so our estimate of
$Z_{2p+1,2p+1}$) will also be small. Empirically we also find that the
resulting matrices are close to the $n \times n$ solution of the CTM equations.

Ideally, we would find approximate solutions at size $n \times n$ by
using (\ref{eq:expand}) to construct very large transfer
matrices, then transform and truncate these transfer matrices according to
(\ref{eq:reduce}). However, since the matrices double in size for each use of
(\ref{eq:expand}), this quickly becomes impractical.

Instead, we transform and truncate the matrices every time we expand. Since we
may choose exactly how many eigenvalues to discard, we can take matrices of
size $n\times n$, expand them to size $2n\times 2n$ and then truncate back down
to size $n \times n$. In this way we can produce a sequence of $n \times n$
matrices which appear to converge to the finite-size solution of the CTM
equations.

As far as we know, it is not proved that this procedure will converge at all,
let alone to the finite-size solution. However, this happens empirically, and as
stated above such a proof is not necessary to establish a rigorous lower bound
--- we can simply observe that the resulting $F$ matrices result in good bounds.
In full, the CTMRG method is as follows:
\begin{enumerate}[(1)]
\item Start with initial $1 \times 1$ matrices $A(a), F(a,b)$ and let $n = 1$.
\item Expand $A$ and $F$ using \Ref{eq:expand}.\label{step:expand}
\item If required, increase $n$ by 1.
\item Diagonalise $A_l(a)$. Let $P(a)$ be the matrix consisting of the
eigenvectors corresponding to the $n$ largest eigenvalues of $A(a)$.
\item Reduce $A_l$ and $F_l$ using \Ref{eq:reduce}.
\item Normalise $A$ and $F$ so that the top-left elements of $A(0)$ and $F(0,0)$
are both 1.
\item Go to step~\ref{step:expand}.
\end{enumerate}
Once the above procedure has converged sufficiently at a given size $n$, we can
substitute the resulting $F$ into the CTM equations to obtain $\xi$ and
$\eta$ and the lower bound. In fact, we can terminate the procedure at any
stage, since any $F$ will result in a bound, but in practice we wait until the
procedure has converged to a desired precision.

We calculate our estimates of $\kappa$ by the formula \cite{CTM:1}

\begin{align}
\label{eqn kappa est}
\kappa &= 
\frac{\left(\sum_a \Tr A^4(a)\right) \left(\sum_{a,b,c,d} \Tr \omega \begin{pmatrix}a&b\\c&d\end{pmatrix}A(a) F(a,c) A(c) F(c,d) A(d) F(d,b) A(b) F(b,a)\right)}{\left(\sum_{a,b} \Tr A^2(a) F(a,b) A^2(b) F(b,a)\right)^2}. 
\end{align}

Each term in this formula represents the partition function of the entire plane.
However, the term in the denominator contains an extra row compared to the first
term in the numerator, while the second term in the numerator contains both an
extra row and column. The net effect is to cancel out everything but the effect
of a single cell. In practice, using this formula is faster and more accurate
than the lower bound computation, but of course does not provide a bound.

\subsection{Asymmetry and the CTMRG}

\begin{figure}
 \begin{center}
  \includegraphics[height=4cm]{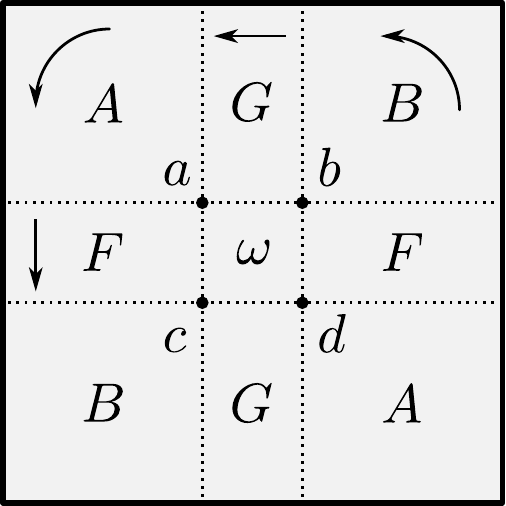}
 \end{center}
  \caption{For a model with symmetry under rotation by $\pi$ we need two sets
of corner transfer matrices $A$ and $B$ as well as half-row and half-column
transfer matrices $F$ and $G$.}
  \label{fig sweep asym}
\end{figure}

The decomposition in (\ref{eq:ctm3}), and thus the CTMRG, implicitly assumes that the corner transfer matrices are equal under a rotation of
$\pi/2$. This naturally occurs if the model itself is invariant under this rotation, which holds true for the majority of the models we analyse. The
exceptions are the RWIM and $Q$-charge models, which are only symmetric
under rotation by $\pi$. To accommodate this, we define
separate half-row and half-column transfer matrices, and two different families
of corner transfer matrices --- see Fig.~\ref{fig sweep asym}.

We will use $A$ to denote the corner transfer matrices for the second
and fourth quadrants of the lattice, and $B$ for the first and third. Similarly
we 
let $F$ and $G$ denote the half-row and half-column transfer matrices respectively. We expand
these matrices in a similar way to \Ref{eq:expand}, taking care with the
orientations:
\begin{subequations}
\begin{align}
 A_l(d)|_{b,c} &= \sum_{a} \omega \abcd G(b,a) \,A(a)\, F(a,c) 
 & F_l(c,a)|_{d,b} &= \omega \abcd F(d,b) \\
 B_l(c)|_{d,a} &= \sum_{b} \omega \abcd F(d,b) \,B(b)\, G(b,a) 
 & G_l(d,c)|_{b,a} &= \omega \abcd G(b,a).
\end{align}
\end{subequations}
To reduce the matrices, we observe that we can represent the transfer matrix of
the entire plane as
either $(AB)^2$ or $(BA)^2$, and these are not necessarily equal. In order to
diagonalise both these matrices and simultaneously preserve the CTM equations,
we define $P_{AB}(a)$ and $P_{BA}(a)$ to be the matrices consisting of the first
$n$ eigenvectors of $A_l(a)B_l(a)$ and $B_l(a)A_l(a)$ respectively. Our reducing
transformations then become
\begin{subequations}
\begin{align}
A(a) &\mapsto P_{AB}^T(a) \,A_l(a)\, P_{BA}(a) 
& F(a,b) &\mapsto P_{BA}^T(a) \,F_l(a,b)\, P_{BA}(b) \\
B(a) &\mapsto P_{BA}^T(a) \,B_l(a)\, P_{AB}(a) 
& G(a,b) &\mapsto P_{AB}^T(a) \,G_l(a,b)\, P_{AB}(b). 
\end{align}
\end{subequations}
The remainder of the method is unchanged.

\subsection{States on bonds}
The CTMRG is defined on a model where the states are located on the vertices of
the $\mathbb{Z}^2$ lattice. As this lattice is self-dual, it can be applied
without change to a model where the states are on the faces. However, we must
make some small adjustments when the states lie on the bonds, as for the $Q$-charge
model. As this model is not invariant under a rotation of $\pi/2$, we will
also use the machinery developed in the previous section.

In order to apply the CTMRG to this form of model, we use a form of the method
similar to that developed by Foster and Pinettes \cite{CTM:RG-FP1,
CTM:RG-FP2}. The face weight $\omega$ becomes a weight around a vertex, and we
shift all the matrices by half a length unit (see Fig. \ref{fig
corner bond}). This also has the effect that the corner transfer matrices are no
longer parameterised by a state --- so we only have one $A$ and one $B$.
Similarly the half-row and half-column matrices, $F(a)$ and $G(a)$, are now
parameterised by a single state.

\begin{figure}
 \begin{center}
  \includegraphics[height=4cm]{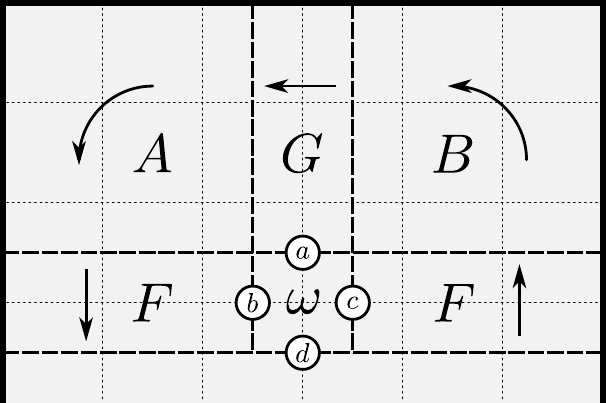}
 \end{center}
  \caption{For a model with bond states we define the corner transfer
matrices $A$ and $B$, and the half-row and half-column matrices $F$ and $G$ as
shown. The edges of the lattice are shown in thin dotted lines. The bond-states surrounding the weight $\omega$
are shown as labelled white circles.}
  \label{fig corner bond}
\end{figure}

We make obvious adjustments to the expansion and contraction of the matrices. These now become
\begin{subequations}
\begin{align}
A_l|_{c,d} &= \sum_{a,b} \omega \abcdvert G(a) \,A\, F(b) 
& F_l(b)|_{d,a} &= \sum_c \omega \abcdvert F(c) \\
B_l|_{d,b} &= \sum_{a,c} \omega \abcdvert F(c) \,B\, G(a) 
& G_l(d)|_{c,b} &= \sum_a \omega \abcdvert G(a) \\
A &\mapsto P_{AB}^T \,A_l\, P_{BA} 
& F(a) &\mapsto P_{BA}^T \,F_l(a)\, P_{BA} \\
B &\mapsto P_{BA}^T \,B_l\, P_{AB} 
& G(a) &\mapsto P_{AB}^T \,G_l(a)\, P_{AB} 
\end{align}
\end{subequations}
where the face weight $\omega$ is given by (\ref{eq:charge_weight}).

Note that the corner transfer matrix formalism which leads to the rigorous lower
bound only requires that the original column transfer matrix $V$ is
symmetric. In the case of the charge model, the symmetry of $V$ follows
from the fact that if you invert all the vertex spins in a column, then the column is
still legal. However, it is no longer the case that $F^T(a)=F(a)$ and our
eigenvalue equations must be changed to reflect this:
\begin{subequations}
\begin{align}
  \xi X &= \sum_a F(a)\, X \,F^T(a) \\
  \eta Y(a) &= \sum_{b,c,d} \omega\abcdvert
  F(b) \,Y(d)\, F^T(c).
\end{align}
\end{subequations}
In practice, applying the power method using these equations can be
unstable for the $Q$-charge model for reasons we have not been able to determine.
In any case, we are able to obtain significantly better bounds for this
model through a relation to the even model (see Section~\ref{sec charge even}).

\section{Results}\label{sec:results}

The CTMRG algorithm described in the previous section was implemented in
\texttt{C++} using the \texttt{Eigen} numerical linear algebra
library \cite{eigenweb}. The main advantage of \texttt{Eigen} over other
libraries is that it readily supports multiple precision computations through
the MPFR library \cite{MPFR}\footnote{At the time of writing the Eigen and MPFR
libraries were available from \texttt{http://eigen.tuxfamily.org} and
\texttt{http://www.mpfr.org/}.}. The eigenvalues of the corner transfer matrices
range over many orders of magnitude, so all computations must be
done at very high precision. For example, for the hard squares model, the
eigenvalues of $A(0)$ and $A(1)$ at size $50\times 50$ range from $2^0$ down to
approximately $2^{-50}$.

All of the computations were run on modest desktop computers running Linux and
took between a few hours and a few days to complete. The fastest computation
(overall) was hard squares and the slowest was charge(3). Indeed, our most
precise results were for hard squares and the least precise were for charge(3). 

Matrix diagonalisation requires $O(n^3)$ operations, and so a full
iteration of the algorithm described in Section~\ref{sec ctmrg} also requires
$O(n^3)$ operations. We have observed that for most of the models, the number of
iterations required to converge at size $n$ does not increase with $n$. However,
since our computer code waits until convergence at every size before proceeding,
we observe that the algorithm is roughly $O(n^4)$.
Note however, that this is the time taken to converge at a fixed matrix size, not to
a fixed number of digits in $\kappa$.

The rate of convergence of our estimates and bounds depends upon the
eigenspectrum of the corner transfer matrices. For the purposes of this
discussion we will assume that the matrices at each size have converged to their
limiting value under the renormalisation procedure. A key observation is that the largest
eigenvalues do not change much with matrix size --- the 4 eigenvalues of $A(0)$
at size 4 are nearly the same as the 4 dominant eigenvalues of $A(0)$ at size
40, and presumably at size 400.

Thus, what changes as we increase matrix size are the eigenvalues we discard by
truncating the matrices. Since each term in our computation of $\kappa$ (via \Ref{eqn kappa est}) is the
trace of the product of 4 corner transfer matrices (some interlaced with
half-row and half-column transfer matrices), we expect that the contribution
of each eigenvalue to the estimate of $\kappa$ is roughly proportional to
its fourth power.

The question of the tightness of the bound is a little more complicated, as it
depends on how close our $F$ matrices are to the solution of the CTM equations.
If they exactly solve the CTM equations at finite size (for some values of the
other matrices), then there is no difference between our bound calculation
using \Ref{kappa bound} and our estimate using \Ref{eqn kappa est}, and so the
error should be the same. As we run the algorithm to convergence at each finite
size, we expect that the extra error induced by having a not-quite-optimal $F$
is minimal.

This analysis shows that the utility of the method hinges on the distribution of
the eigenvalues of the corner transfer matrices. If the eigenvalues decay
slowly, then by truncating the matrices we will throw away quite large
contributions to the estimates and bounds of $\kappa$. In such a situation, we would expect
our estimates and bounds to converge slowly with increasing matrix size. On the other hand,
if the eigenvalues decay quickly then the truncated eigenvalues will have
little impact.

It is believed \cite{CTM:RG-Eigenvalues} that the eigenvalues of the corner
transfer matrices scale as
\begin{equation}
\lambda_k \sim \exp \left(-c (\ln k)^2 \right) 
\label{eq:eigvals}.
\end{equation}
Here we have denoted the $k^\text{th}$ largest eigenvalue by $\lambda_k$; since these
do not vary significantly between different matrix sizes, we have not specified
the dimension of the underlying matrix.

\begin{figure}
\centering
\includegraphics[height=5cm]{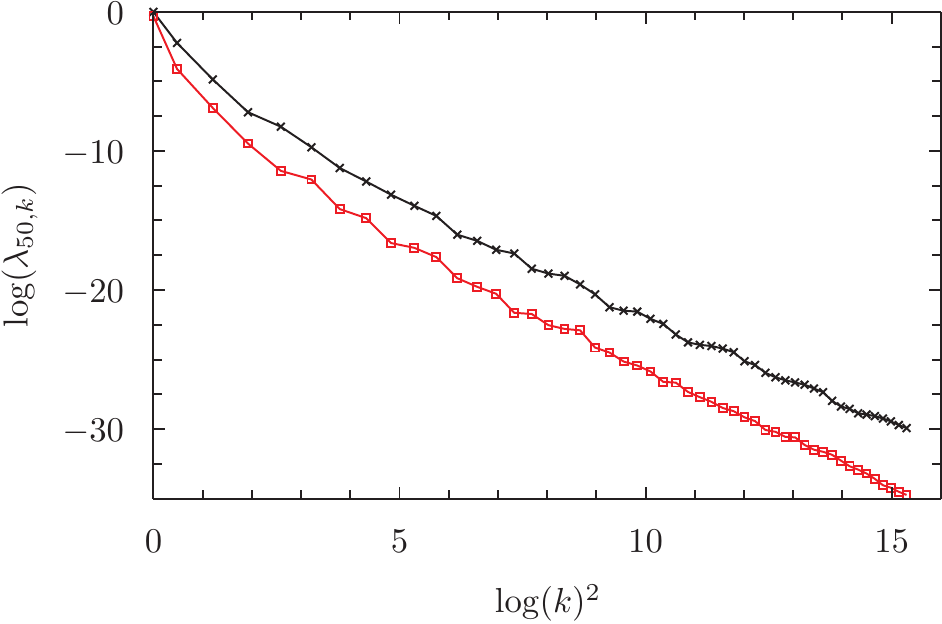}
\caption{A plot of $\log(\lambda_k)$ vs $\log(k)^2$, for the eigenvalues of
$A(0)$ (upper curve) and $A(1)$ (lower curve) at size 50 for the hard squares model. For
larger $k$, the plots appear approximately linear, supporting the conjectured
scaling form of \Ref{eq:eigvals}.} 
\label{fig ev spec}
\end{figure}

If the above scaling form holds true --- and our observations suggest that it
does (see Fig.~\ref{fig ev spec}) --- then we expect the error in our
estimates to have a similar asymptotic form with increasing matrix size.  Thus,
to achieve an error approximately $\varepsilon$, we need to take matrices of
size $k\times k$ with 
\begin{align}
  k & \sim \exp\left( \alpha \sqrt{ - \log( \varepsilon )} \right).
\end{align}
Substituting $\varepsilon = 2^{-N}$ gives
\begin{align}
  k & \sim \exp\left( \beta \sqrt{ N } \right).
\end{align}
Since the running time and memory required grow polynomially with matrix size,
the above argument supports the assertion that the time and memory
required by the algorithm to correctly compute $N$ digits grows
subexponentially with $N$.

\begin{table}
\renewcommand{\arraystretch}{1.1}
\begin{center}
\caption{Lower bounds and estimates for the growth rates of the various models.
For each model we give the full decimal expansion of the lower bound, and then we give the last few digits of our
estimate in brackets. These replace that many digits in the lower bounds --- the preceding digits are identical. We underline the
digits which agree with the previous best known lower bounds
(references given underneath the model name).}
\begin{tabular}{|c||c|l|}
\hline
Model & Matrix size & Lower bound on (and estimate of) $\kappa$ \\
\hline
Hard squares & 256 & \underline{1.503 048 082 475 332 264 322 0}66 329 475 \\
\cite{Friedland2010} & & \hspace{1ex} 553 689 385 781 038 610 305 062 028 101 \\
& & \hspace{1ex} 735 933 850 396 923 440 380 463 299 47 (65)\\
\hline
NAK & 256 & \underline{1.342 643 951} 124 601 297 851 730 161 875\\
\cite{louidor2010} & & \hspace{1ex} 740 395 719 438 196 938 393 943 434 885\\
& & \hspace{1ex} 455 0 (1)\\
\hline
RWIM & 128 & \underline{1.448 957 3}71 775 608 489 872 231 406 108\\
\cite{louidor2010} & & \hspace{1ex} 136 686 434 371 (7)\\
\hline
\hline
Even & 128 & \underline{1.35}7 587 502 184 123 (5) \\
\cite{louidor2010} & & \\
\hline
Charge$(3)$ & 74 & (\underline{1.3}57 587 50) \\
\cite{louidor2010} & & \\
\hline
\hline
4-Colouring & 96 & \underline{2.336 056 64}1 041 133 656 814 01 (4)\\
\cite{Lundow} & & \\
\hline
5-Colouring & 64 & \underline{3.250 404 923 16}7 119 143 819 73 (6)\\
\cite{Lundow} & &\\
\hline
\end{tabular}
\end{center}
\label{tab:kappa}
\end{table}

In Table~\ref{tab:kappa}, we list our lower bounds and estimates for the growth rates of the
models we have considered. Due to the instability of the bound calculation for the charge(3) model, we only give the estimate for this model. Note that since RWIM is
asymmetric, we computed the bounds and estimates both using
\Ref{eq:ctm1} and~\Ref{eq:ctm2} as stated, and also with $F$ replaced
by $G$; they gave the same bounds (to the precision stated in the table).

To show how much our bounds are an improvement on existing bounds, we have underlined the digits which agree with previous best known bounds (after which, our bounds are greater). Of course, rigorously measuring the accuracy of a bound is impossible without exact knowledge of the quantity being estimated, or at least a comparable upper bound. However, comparing the number of digits with which our estimates agree with our bounds and previously known bounds, we are confident that our bounds are a great improvement in all cases.

We note that our estimates agree with all 43 digits of the estimate for the hard 
squares growth rate given in \cite{BaxterPlanar}, and agree up to our last 2 
digits with a naive evaluation of the $q$-colouring series given in 
\cite{Jacobsen} with $q = 4,5$. We have attempted to identify all of these 
constants using the Inverse Symbolic Calculator \cite{iscweb}. Unfortunately, we 
were not able to identify any of them.

Our estimates lead to the rather surprising observation that the growth rates of
the charge(3) and even models are identical --- at least to the 8 significant
digits we have computed. This observation led us to Conjecture~\ref{conj charge
even}. We discuss this further below.

\section{The charge$(3)$ and even models}\label{sec:charge-even}
\label{sec charge even}

As noted above, our estimates of the growth rates of the charge$(3)$ and even
models agree to 8 significant digits, which is the limit of our accuracy
for the charge$(3)$ model. Because of this we conjecture that they are equal.
Should this conjecture be true, our lower bound and estimate for the charge$(3)$ model can be
extended by a few digits, using values from the even model. 

\begin{figure}
 \centering
  \includegraphics[height=3cm]{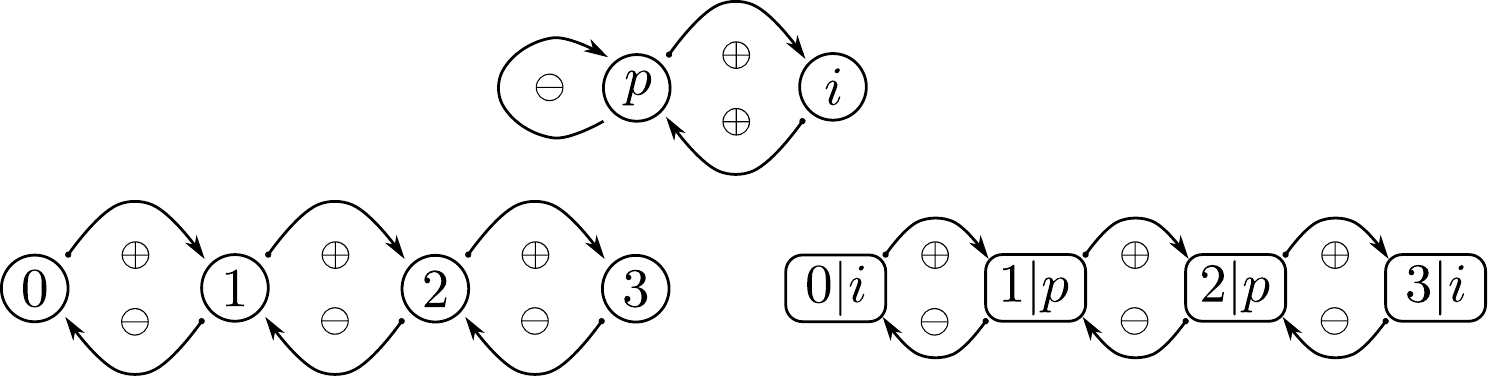}
  \caption{The allowed transitions for the one-dimensional \pdi (top) and
$Q$-charge (left) models. There is a 2-1 graph homomorphism from the
$Q$-charge finite-state graph to that of the \pdi
model --- namely map the nodes $0,3$ to the node $i$ and the nodes $1,2$ to
the node $p$.}
  \label{fig qpi map}
\end{figure}

The one-dimensional counterpart to this conjecture can be proved quite readily:
\begin{theorem}
For the one-dimensional lattice,
\begin{align}
 \kappa_{charge(3)} &= \kappa_{even} = \frac{1+\sqrt{5}}{2}.
\end{align}
\end{theorem}
\begin{IEEEproof}
Perhaps the most obvious way to demonstrate this is to construct the
appropriate transfer matrices for the $Q$-charge and \pdi models (see
Fig.~\ref{fig qpi map}). These are
\begin{align}
  V_{\mbox{\scriptsize \pdi}} &= \begin{pmatrix}
             1 & 1 \\
	     1 & 0
            \end{pmatrix}
 & V_{\mbox{\scriptsize $Q$-charge}} &= \begin{pmatrix}
		  0 & 1 & 0 & 0 \\
		  1 & 0 & 1 & 0 \\
		  0 & 1 & 0 & 1 \\
		  0 & 0 & 1 & 0 
		  \end{pmatrix}.
\end{align}
A quick computation shows that the dominant eigenvalue (and hence the growth rate)
of each is $\frac{1+\sqrt{5}}{2}$.

It is more instructive, however, to construct a mapping between the two
states-on-bonds models. This mapping is a $2$-$1$ graph homomorphism between the
two underlying state graphs as illustrated in Fig.~\ref{fig qpi map}. Simply
map the vertices labelled by $0,3$ to the vertex labelled $i$, and the vertices
labelled $1,2$ to the vertex labelled $p$.

Now consider any valid $Q$-charge configuration in one dimension. Using the
above mapping we can transform it into a valid \pdi configuration. The reverse
mapping is $1$ to $2$. To see this pick a single $p$ in the \pdi configuration
--- we may choose whether this maps to a $1$ or a $2$ in the corresponding
$Q$-charge configuration. Once this choice is made, the rest of the mapping is
fixed. The result follows.
\end{IEEEproof}

The mapping used in the above proof can also be used to show that the growth rate
of the charge(3) model in two dimensions is at least that of the even model.
\begin{theorem}
For the two-dimensional lattice,
\begin{align}
	\kappa_{charge(3)} &\geq \kappa_{even} (\geq 1.357\,587\,502\,184\,123).
\end{align} 
\end{theorem}
\begin{IEEEproof}
Choose any bond on the lattice, and draw a line parallel to $y = x$ through this bond. Choose all bonds intersecting this line to have $Q$-charge states 0 or 2, as shown in Fig. \ref{fig:parity}. Now consider a valid $Q$-charge configuration conforming to these restrictions. The set of such configurations is a (strict) subset of all valid $Q$-charge configurations. However, since the bond-states must increase or decrease by 1 as we read across horizontal and vertical lines, we have determined the parity of \emph{every} bond state on the lattice.
\begin{figure}\centering
\includegraphics[scale=0.7]{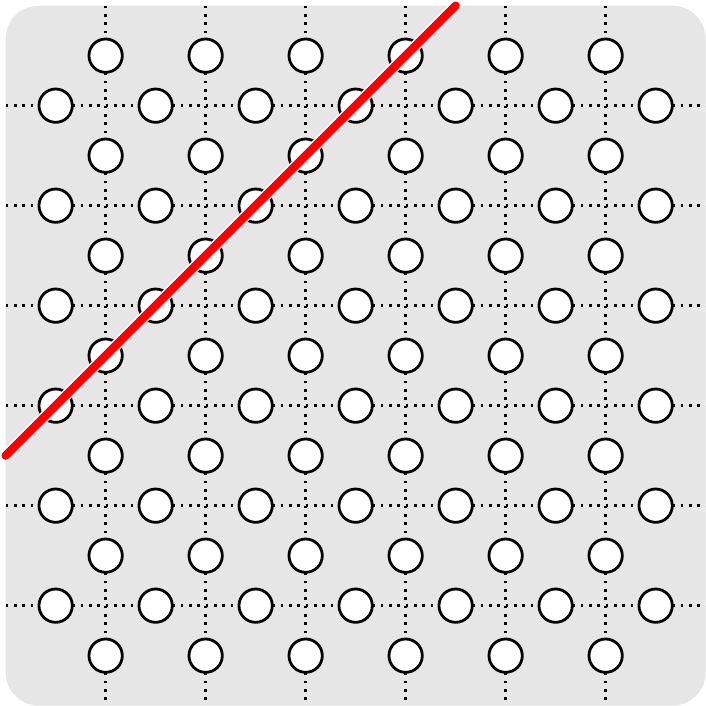}
\caption{Restriction for the $Q$-charge $\leftrightarrow$ \pdi mapping. All
states on the red line must have state 0 or 2 in the $Q$-charge model.}
\label{fig:parity}
\end{figure}

Now in every cell, the parities of the N and W bonds are determined and identical. A full enumeration of the (five) possible cell configurations shows that in this case, the above mapping between the $Q$-charge and \pdi models is a bijection between all valid $Q$-charge and \pdi cell configurations. Thus the number of configurations for each model is the same. Since we have only taken a subset of valid $Q$-charge configurations, but all valid \pdi configurations, this provides the inequality.
\end{IEEEproof}

Unfortunately, we have not been able to prove the reverse inequality, though we can make a non-rigorous argument as to why this should be so. In every $Q$-charge configuration, there will be a number of faces where the N and W bonds have identical parity (the `good' faces), and a number of faces where the N and E bonds have the same parity (the `bad' faces). As mentioned above, it is possible to take configurations which consist entirely of good faces, and these are equivalent to the \pdi configurations.

It is also possible to take configurations which consist entirely of bad faces, and a similar analysis of cell configurations shows that these faces map onto (invalid) \pdi faces containing exactly one $i$ and three $p$'s. By considering an $i$ as an `occupied' bond and a $p$ as an `vacant' bond, we see that these configurations are in fact fully-packed dimer configurations. The growth rate of such configurations has been proved \cite{Kasteleyn} to be
\begin{align}
\kappa & = e^{Catalan/\pi} = 1.338\,515\,15\dots
\end{align}
which is less than our proved lower bound for $\kappa_{even}$. Therefore, these configurations are exponentially infrequent in the set of $Q$-charge configurations.

Another possibility is to take configurations where the good faces and bad
faces are arranged in a checkerboard pattern. We call this model the `elbow'
model, due to the characteristic shape of the individual bond configurations. The growth rate of this model can be rigorously bounded above using the
method of Calkin and Wilf \cite{CalkinWilf}. Using transfer matrices of width $8$, we calculate
$\kappa_{elbow} \leq 1.333\,521\ldots$, which is again less than the proved
lower bound for $\kappa_{even}$. 

\begin{figure}
\centering
\includegraphics[height=8cm]{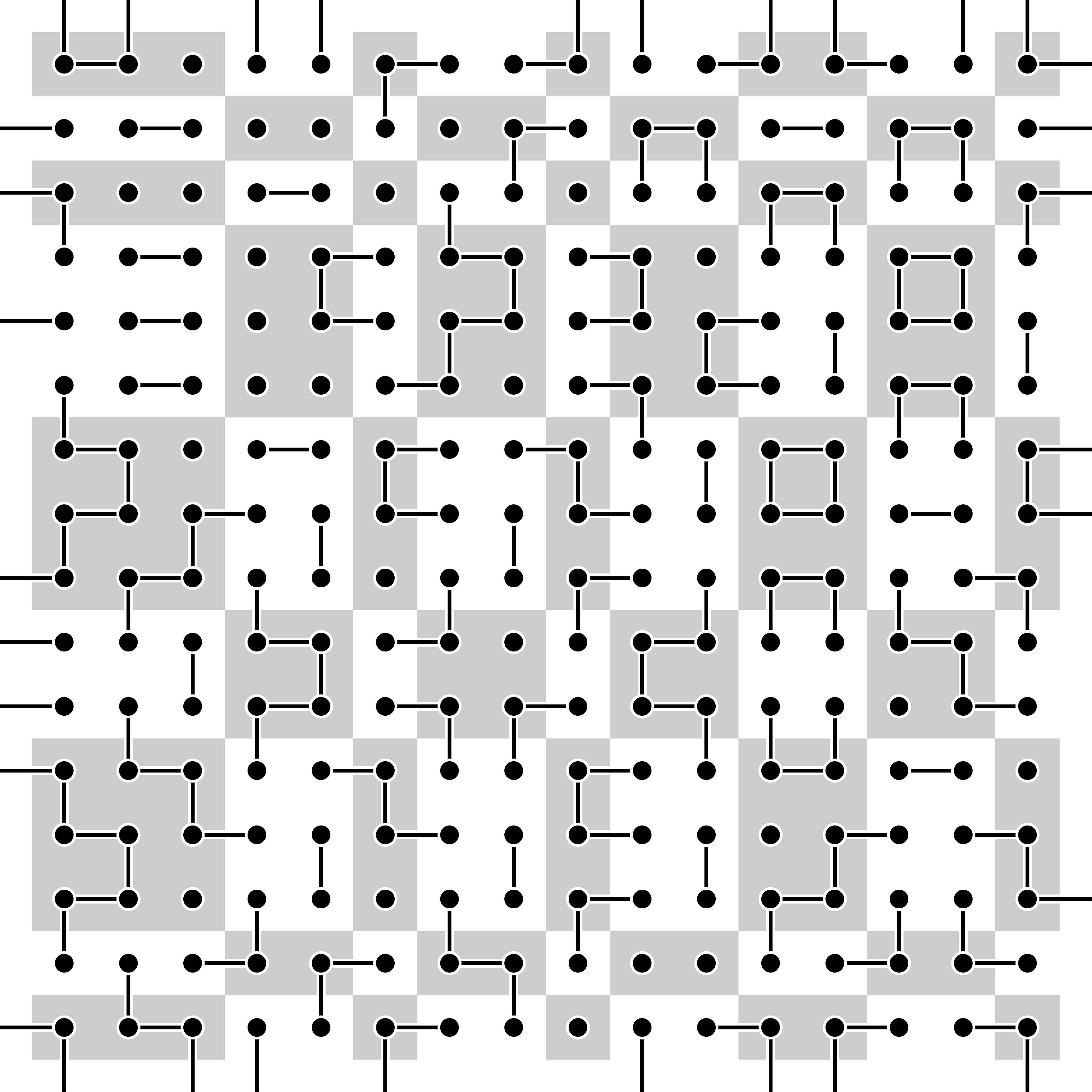}
\caption{An example of the $Q$-charge model mapped to a
`tartan' combination of \pdi and fully-packed dimer constraints. The shaded regions obey
the \pdi constraint, while the unshaded regions obey the dimer constraint.}
\label{fig:tartan}
\end{figure}

Now choose any fixed set of parities for the rows and columns of the lattice.
In these configurations, the good faces and bad faces form a `tartan'-like
configuration
(see
Fig.~\ref{fig:tartan}). The growth rate of such configurations is some combination
of the growth rates of the even model (the good faces), the fully-packed dimer
model (the bad faces) and the elbow model (an interfacial term
between good and bad faces). Out of these three models, the even model has the
largest growth rate and so the growth rate of such configurations cannot be greater than
$\kappa_{even}$.

Now there are $2^{2\sqrt{N}}$ ways to choose these parities on a lattice of $N$ spins, a sub-exponential factor. Each of these ways results in a growth rate not more than $\kappa_{even}$, and therefore the set of all valid $Q$-charge configurations has a growth rate not greater than $\kappa_{even}$.

This (admittedly non-rigorous) argument, combined with the accuracy of our numerical results and the above theorem, presents what we find to be a compelling case that $\kappa_{charge(3)} = \kappa_{even}$.

\section{Conclusion}\label{sec:conclusion}

In this paper, we have constructed very accurate rigorous lower bounds on the
growth rate of several families of data storage models: exclusion, colouring,
even, and charge(3) models. The process we use also gives us very precise
estimates of the growth rates. This is achieved using Baxter's corner transfer
matrix formalism combined with the corner transfer matrix renormalisation group
method of Nishino and Okunishi.

The precision of our bounds allows us to conjecture the equality of the growth rate of two of our models, the charge(3) and even models. We have proved one inequality for this relation, and have a non-rigorous (but overall convincing) argument for equality.

For future work, we would like to be able to construct upper bounds of similar quality. The quality of proved lower bounds have typically been much higher than that of upper bounds. It may be possible to use linear algebra results to obtain a bound on the error of our estimates, and so derive an upper bound for the growth rate. It certainly is possible to bound the distance between our estimate and the nearest eigenvalue of the transfer matrix using a theorem of Wilkinson \cite{Wilkinson}. However, there is no guarantee that this eigenvalue is the dominant eigenvalue, and so more work needs to be done to establish a provable upper bound.

In addition, there are still a number of constraints which we have not analysed in this paper, for example higher-order charge or colouring models. These are amenable to identical analysis, and are the subject of future work.

\section*{Acknowledgements}

Much of this work took place while YBC was a postdoctoral fellow at the 
University of Vienna. YBC acknowledges financial support from the University of 
Vienna and thanks the University of British Columbia for their hospitality. AR 
acknowledges financial support from NSERC and thanks the University of Vienna 
for their hospitality during his visit. Both authors thank Brian Marcus for many 
helpful discussions and proof-reading.

\bibliographystyle{IEEEtran}
\bibliography{lb_refs}



\end{document}